\documentstyle{article}
\renewcommand{\baselinestretch}{1.0}
\title{On framed instanton bundles and their deformations}
\author{Andreas Matuschke\\
        \small Humboldt-Universit\"{a}t zu Berlin,
         Graduiertenkolleg ``Geometrie und Nichtlineare Analysis''\\
         \small Ziegelstrasse 13 a, D-10099 Berlin\\
         \small email: matuschk@mathematik.hu-berlin.de}
          
\date{October 10, 1996}

\newcommand{\ZZ}{\rm Z\!\! Z}

\newcommand{\RR}{\rm I\! R}

\newcommand{\CC}{\rm C\!\!\! I\,}
\newcommand{\PP}{\rm I\! P}
\newcommand{\EE}{\rm I\! E}

\newcommand{\Hom}{{\cal H}\hspace{-0.1em}{\it o\hspace{-0.1em}m}}
\newcommand{\End}{{\cal E}\hspace{-0.1em}{\it n\hspace{-0.1em}d}}

\newtheorem{theorem}{Theorem}[section]
\newtheorem{lemma}[theorem]{Lemma}
\newtheorem{proposition}[theorem]{Proposition}
\newtheorem{corollary}[theorem]{Corollary}

\newcommand{\ssty}{\scriptstyle}

\begin{document}

\maketitle

\renewcommand{\baselinestretch}{1.0}
\small \normalsize
\begin{abstract}
\noindent  
We consider a compact twistor space P and assume that
there is a surface $\rm S\subset P$, which has degree one
on twistor fibres and contains a twistor fibre F, e.g.
P a LeBrun twistor space (\cite{lebrun},\cite{kurke1}).  
Similar to \cite{donaldson} and \cite{buchdahl2} we examine
the restriction of an instanton bundle V equipped with
a fixed trivialisation along F to a framed vector bundle over (S,F).
First we develope inspired by \cite{huy-lehn1} a suitable 
deformation theory for vector bundles over an 
analytic space framed by a vector bundle over a subspace
of arbitrary codimension. 
In the second section we describe the restriction as
a smooth natural transformation into a fine moduli space. 
By considering framed $\rm U(r)-$instanton bundles 
as a real structure on framed instanton bundles over P, we show that
the bijection between isomorphism classes of framed 
$\rm U(r)-$instanton bundles
and isomorphism classes of framed vector bundles over (S,F) due to 
\cite{buchdahl2} is actually an isomorphism of moduli spaces.
\indent
\end{abstract}

\renewcommand{\baselinestretch}{1.5}
\small \normalsize

\section{Introduction}

We consider the twistor fibration 
$\rm \pi : P \rightarrow M$
over a real four-dimensional compact manifold M 
with self-dual Riemannian metric. 
P is a three-dimensional complex manifold with an induced antiholomorphic fixpoint free involution $\rm \tau$, 
an antipodal map on the twistor fibers 
(cf. \cite{AHS}, \cite{buchdahl1}, \cite{friedrich}).
A line on P is a complex submanifold 
$\rm L\subset P$ with $\rm L\cong \PP_{\CC}^1$ and normal bundle
$\rm {\cal N}_{L|P} \cong {\cal O}_{\PP^1}( 1 ) ^{\oplus 2}$. 
In particular, twistor fibres are lines. We denote with
$\rm \mu : Z \rightarrow H$ the universal Douady-family of lines in P. 
The involution $\rm \tau$ maps lines to lines and consequently induces 
an antiholomorphic involution on H. 
Then M appears as a set of fixpoints of $\rm \tau$ 
and moreover as a real-analytic submanifold of H (cf. \cite{AHS}, \cite{kurke2}):
\[
\begin{array}{ccccc}
&&\rm H \times P&&\\
&&\cup&&\\
\rm P=Z {\times}_H M& \longrightarrow &
\rm Z& \stackrel{\nu}{\longrightarrow}&
\rm P\\
\rm\big\downarrow {}^\pi&&\big\downarrow {}^\mu&&\\
\rm M&
\rm  \subset &
\rm H &&
\end{array}
\]

Let $\rm S \subset P$ be a surface intersecting twistor fibres with multiplicity 1 and containing a twistor fibre F. 
By \cite{kurke1}, Prop.2.1, S is a smooth algebraic surface and F the only twistor fibre in S. 
With $\rm \bar{S} = \tau \left( S \right)$ we have $\rm F=S \cdot \bar{S}$.
The linear system $\rm |F|$ defines a birational morphism $\rm S\rightarrow \PP_{\CC}^2$ and we have  
$\rm M\approx \,\stackrel{n}{\#} \left( -\PP^2_{\CC} \right)${} or $\rm S^4$.
A well examined class of examples are the LeBrun twistor spaces,
which fulfill the additional property $\rm dim|S|\ge 1$, cf. \cite{lebrun},
\cite{kurke1}. 
These twistor spaces are classified as modifications of conic bundles and 
are algebraic in the sense of M. Artin \cite{knutson}.

A (mathematical) instanton bundle is an holomorphic vector bundle on P trivial
on all twistor fibres. 
The Penrose-Ward transformation gives us an analytic equivalence between the categories of instanton bundles and of pairs 
$\rm \left( E, \nabla \right)$ of complex vector bundles on M with self-dual connection.
The pair $\rm \left( E, \nabla \right)$ is associated to the 
$\rm C^{\infty}$-bundle $\rm \pi^{\ast}E$ together with the holomorphic structure defined by $\rm \bar{\partial}=\left( \pi^{\ast}\nabla \right)^{01}$.
Conversely, an instanton bundle V gives rise to a pair 
$\rm \left( E, \nabla \right)$ by taking 
$\rm E=\left( \mu_{\ast}\nu^{\ast}V \right) \big{|}_M$ and 
$\rm \nabla$ as restriction of
\[
\begin{array}{ccc}
\rm \mu_{\ast}\left( {\cal O}_Z \otimes_{\nu^{-1}{\cal O}_P} \nu^{-1}V \right)&
\stackrel{\rm \mu_{\ast}( d_{Z|P}\otimes id_{\nu^{-1}V} )}{-\!\!\!-\!\!\!-\!\!\!-\!\!\!-\!\!\!-\!\!\!-\!\!\!-\!\!\!- \!\!\!-\!\!\!-\!\!\!-\!\!\!-\!\!\!-\!\!\!-\!\!\!-\!\!\!\rightarrow}&
\rm \mu_{\ast}( {\Omega}_{Z|P}^1 \otimes_{{\cal O}_Z} \nu^{\ast}V )\\
\Big\downarrow \nabla&
&
\Big\downarrow \cong\\
\rm {\Omega}_H^1 \otimes_{{\cal O}_H} \mu_{\ast}\nu^{\ast}V &
\stackrel{\cong}{-\!\!\!-\!\!\!- \!\!\!-\!\!\!-\!\!\!-\!\!\!-\!\!\!-\!\!\!-\!\! \!-\!\!\!-\!\!\!-\!\!\!-\!\!\!-\!\!\!-\!\!\!-\!\!\!\rightarrow}&
\rm \mu_{\ast}{\Omega}_{Z|P}^1 \otimes_{{\cal O}_H} \mu_{\ast}\nu^{\ast}V 
\end{array}
\]
(cf. \cite{AHS}, \cite{kurke2}).

For G a linear group as for example $\rm U(r)$, 
$\rm SU(r)$, $\rm Sp(r)$ or $\rm SO(r)$, 
a G-instanton on M is given by a complex G-vector bundle E with self-dual connection $\rm \nabla$ compatible with the G-structure on E. 
The Penrose-Ward transformation associates a G-instanton to an instanton bundle with additional properties. 
In particular, the category of $\rm U(r)$-instantons
is analytically equivalent to the category of instanton bundles V on P, for which there is an isomorphism 
$\rm \varphi : V {\cong} \tau^\ast \bar{V}^\vee$ with
$\rm \tau^\ast \bar{\varphi}^\vee = \varphi$, where $\rm \bar{V}^\vee$ denotes the bundle of antilinear forms. 
We denote these instanton bundles as 
$\rm U(r)$-instanton bundles or as physical instanton bundles 
(cf. \cite{AHS}, \cite{buchdahl1}, \cite{kurke2}).

For an instanton bundle V on P we can fix a trivialisation 
$\rm \alpha : V|_F \cong {\cal O}_F^r$ along our twistor fibre 
$\rm F=S\cdot\bar{S}$. 
The resulting pair $\rm \left( V,\alpha \right)$ is called a framed instanton bundle. 
The frame lives in codimension 2.
By restricting V to S we obtain a framed vector bundle on a smooth rational surface, framed along a divisor that is big and nef. 
For this case the moduli problem is well examined by Lehn and Huybrechts in the algebraic-projective case (\cite{lehn}, \cite{huy-lehn1}) and by 
L\"ubke \cite{lubke} from the analytic point of view.
It was an idea of Donaldson \cite{donaldson} to use this restriction map
to discuss moduli of framed instantons. 
By the results of Buchdahl \cite{buchdahl2} it follows, that there is a bijection between the isomorphism classes of framed $\rm U(r)$-instantons and framed vector bundles on S.

In the first section this paper we develop the deformation theory 
for framed vector bundles on 
analytic spaces, where the frame has support of arbitrary codimension. 
I.e. we answer questions like: what is the tangent space and 
when is the deformation functor formally smooth. 
Everything in this section can be done in the
same way and with the same results for seperated algebraic spaces 
of finite type over some field of characteristic zero.

In the second section we apply these results to framed instanton 
bundles and describe the mentioned restriction map infinitisimally. 
By describing the functor of framed 
$\rm U(r)$-instanton bundles as a real structure on the functor of framed mathematical instanton bundles, 
we show a real-analytic isomorphism between the moduli of framed 
$\rm U(r)$-instanton bundles and the moduli of framed vector bundles on S.
                                                                                                                                                                                                                                                                                                                                                                                                                                                                                                                                                                                                                                                             
This paper is based on my Diplomarbeit,
which was produced in March 1994
under the supervision of Herbert Kurke in Berlin.

\section{Deformations of framed vector bundles}

\subsection{Families and deformations of framed vector bundles}
The following notation is fixed through the whole section.
Let $\rm X_0$ be an analytic space, 
$\rm Y_0$ a closed subspace with associated ideal sheaf
$\rm {\cal J}_{Y|X}$
and ${\rm W_0}$ a fixed locally free sheaf of mo- dules over $\rm Y_0$ as in  \cite{grauert}. 
Then a {\it framed vector bundle} to the data                        
$\rm (X_0,Y_0,W_0)$  
consists of a pair 
$\rm ( E_0, \alpha_0)$
with ${\rm E_0}$ a locally free module sheaf over ${\rm X_0}$ and
${\rm \alpha_0: E_0|_{Y_0} \rightarrow W_0}$ a framing isomorphism. 
A {\it morphism}  f between two framed vector bundles                                           $\rm  \left(E_0,\alpha_0\right)$ and
$\rm \left(E'_0,\alpha'_0\right) $ is a sheaf morphism 
${\rm f:E_0\rightarrow E'_0}$ with 
${\rm \alpha'_0\circ\left(f|_{Y_0}\right)=\alpha_0}$.
For an analytic space S a {\it family} of framed vector bundles for
 the given data ${\rm \left(X_0,Y_0,W_0\right) }$ parametrized by S is a framed vector bundle 
${\rm \left( E,\alpha \right) }$ for the data 
${\rm \left( S\times X_0, S\times Y_0, p^* W_0 \right) }$, where            
${\rm p:S\times Y_0 \rightarrow Y_0 }$ is the projection.

We fix a framed vector bundle $\rm ( E_0, \alpha_0 )$. 
A {\it deformation} of $\rm ( E_0, \alpha_0 ) $ over 
a germ of an analytic space $\rm (S,s_0)$
is represented by a triple 
${\rm (E,\alpha,\psi ) }$, where
${\rm \left( E,\alpha \right) }$ 
is a family of framed vector bundles over S 
and  ${\rm \psi }$ is an isomorphism from                                  
${\rm E|_{\left\{ s_0 \right\} \times X_0} }$ to 
${\rm E_0 }$, such that the diagram 
\[
\begin{array}{ccc}
{\rm E|_{\left\{ s_0 \right\} \times Y_0 } }&
\begin{array}{c}
{\cong}\\[-3ex]
-\!\!\!-\!\!\!-\!\!\!-\!\!\!\rightarrow\\[-2.5ex]
\psi
\end{array}&
{\rm E_0|_{Y_0}  }\\
{\rm {\cong}{\Big\downarrow}{\alpha} }&
&
{\rm {\cong}{\Big\downarrow }{\alpha_0} }\\
{\rm W_0 }&
\begin{array}{c}
{=}\\[-3ex]
-\!\!\!-\!\!\!-\!\!\!-\!\!\!\rightarrow\\[-2.5ex]
{}
\end{array}&
{\rm W_0 }
\end{array}
\]
commutes. Here, the pointed space $\rm (S,s_0)$ 
is any representative of our
germ. 
If $\rm \eta : (T,t_0)\rightarrow (S,s_0)$ is a
local analytic isomorphism, then 
${\rm ((\eta\times id_{X_0})^{\ast}E, 
(\eta\times id_{X_0})^{\ast}\alpha,(\eta\times id_{X_0})^{\ast}\psi ) }$
represents the same deformation.

Two deformations
$\rm (E,\alpha,\psi )$ and $\rm (E',\alpha',\psi') $ of
$\rm ( E_0, \alpha_0 ) $ over the same germ 
can be realized as families over the same pointed analytic space
$\rm (S,s_0)$.
Then a morphism 
${\rm f\!: (E,\alpha,\psi ) \rightarrow\! (E',\alpha',\psi') }$
is represented by 
a local analytic isomorphism
$\rm \eta :(T,t_0)\rightarrow (S,s_0)$
and a morphism 
\[{\rm f: ((\eta\times id_{X_0})^{\ast}E,(\eta\times id_{X_0})^{\ast}\alpha ) \rightarrow 
((\eta\times id_{X_0})^{\ast}E',(\eta\times id_{X_0})^{\ast}\alpha') }\] 
of framed vector bundles with 
${\rm (\eta\times id_{X_0})^{\ast}\psi' \circ 
( f|_{s_0 \times X_0} ) = (\eta\times id_{X_0})^{\ast}\psi }$.

We denote with ${\rm {\cal M}\left( X_0,Y_0,W_0 \right) }$ 
the functor
${\rm (analytic\; spaces ) \longrightarrow ( sets ) }$
defined by
\[ 
\begin{array}{ccc}
{\rm S}&
\longrightarrow&
\left\{ \begin{array}{c}
{\rm isomorphism\; classes\; of\; families\; of}\\[-1.5ex]
{\rm framed\; vector\; bundles\; to\; the\; data}\\[-1.5ex]
{\rm \left( X_0,Y_0,W_0 \right)\; parametrized\; by\; S} 
\end{array} \right\}
\end{array}
\]
and with ${\rm Def \left( E_0,\alpha_0 \right) }$ the functor 
${\rm \left(germs\; of\; analytic\;spaces\right)\longrightarrow\left(         sets \right) }$\\
defined by
\[ 
\begin{array}{cccc}
{\rm \left( S,s_0 \right)}&
\longrightarrow&
\left\{ \begin{array}{c}
{\rm isomorphism\;classes\;of\;deformations}\\[-1.5ex]
{\rm of\; \left( E_0,\alpha_0 \right)\;over\; \left( S,s_0 \right) } 
\end{array} \right\}&
{\rm .}
\end{array}
\]
For a local ring (R,m) we write often 
${\rm Def \left( E_0,\alpha_0 \right)(R) }$ instead of
${\rm Def \left( E_0,\alpha_0 \right)(SpecR,m) }$.
We note that in the case 
$\rm dim(S,s_0)=0$ the germ represented by $\rm (S,s_0)$ 
coincides with the isomorphism class of this pointed space.

\subsection{Deformation-theoretical tangent spaces}
Let $\rm \CC [ \varepsilon ]$ be the 
$\rm \CC$-algebra 
$\rm \CC [X] / ( X^2)$ with
$\rm \varepsilon = X \,\, mod\, X^2$. 
According to Schlessinger \cite{schlessinger} or M. Artin (\cite{artin1}, \cite{artin2}) we define 
$\rm Def( E_0, \alpha_0 )(\CC[ \varepsilon ])$ 
to be the deformation-theoretical tangent space of 
$\rm {\cal M}( X_0,Y_0,W_0)$
at the point  
$\rm ( E_0, \alpha_0 ) \in {\cal M}( X_0,Y_0,W_0 )( \CC )$. 
$\rm Aut(E_0,\alpha_0)$ acts on 
$\rm Def(E_0,\alpha_0)$ by
$\rm g\cdot (E,\alpha,\psi)=(E,\alpha,g\circ \psi)$.
If this action is trivial and if
there is a local moduli $\cal M$ for 
$\rm ( E_0, \alpha_0 ) $, then we have canonically 
$\rm T_{ ( E_0, \alpha_0 )}{\cal M} = 
Def( E_0, \alpha_0 )( \CC [ \varepsilon ])$.

Due to \cite{schlessinger}, 
$\rm Def ( E_0, \alpha_0 )(\CC \left[ \varepsilon \right])$ 
has a natural structure as complex vector space if the canonical mapping
\[ \rm Def( E_0, \alpha_0 )(\CC [ \varepsilon ]\times_{\CC} \CC
[ \varepsilon ]) \longrightarrow Def ( E_0, \alpha_0)
( \CC [ \varepsilon ])\times Def ( E_0, \alpha_0)(\CC [\varepsilon ])\] 
is bijective.
For $\rm \lambda \in \CC$ and  
$\rm mult(\lambda ) : \CC [ \varepsilon]\rightarrow \CC [ \varepsilon]$ 
defined by 
$\rm x+y\varepsilon \rightarrow x+\lambda y\varepsilon$, 
we obtain
\[ \rm Def( E_0, \alpha_0)(mult ( \lambda )):
Def( E_0, \alpha_0)(\CC [ \varepsilon])\longrightarrow 
Def( E_0, \alpha_0 )( \CC [ \varepsilon ]),\]
which explains the scalar multiplication.
The mapping $\rm add:\, \CC [ \varepsilon ] 
\times_{\CC} \CC [ \varepsilon ]\rightarrow 
\CC [ \varepsilon] $ defined by 
$\rm add( x+y\varepsilon, x+z\varepsilon)= x+( y+z)\varepsilon$ 
induces
\[ \rm Def( E_0, \alpha_0 )(add):
Def( E_0, \alpha_0)(\CC [ \varepsilon] 
 \times_{\CC} \CC[ \varepsilon])
\longrightarrow Def( E_0, \alpha_0)(\CC[ \varepsilon]).\]
Under the bijection above we obtain the additive structure on 
$\rm Def( E_0, \alpha_0)(\CC [ \varepsilon])$. 
The assumed bijection is a consequence of the following lemma.

\begin{lemma}
Let $\rm A'\rightarrow A$ be a small extension of Artin rings. Then for any morphism $\rm B \rightarrow A$ of Artin rings and for $\rm B'=B\times_{A}A'$
the canonical mapping
\[
\rm Def( E_0, \alpha_0)( B') \longrightarrow 
Def( E_0, \alpha_0)(A')\times_{ Def(E_0,\alpha_0)(A)} 
Def( E_0, \alpha_0)(B)\] 
is bijective.
\end{lemma}

{\it Proof:} 
We recall that an {\it Artin ring} is a local $\rm \CC$-algebra of
finite dimension. 
A {\it small extension} 
$\rm A'\rightarrow A$ is given, 
if there is an element $\rm t\in A'$, 
such that $\rm A= A'/tA'$ and 
$\rm t\cdot m_{A'}=0$. 
Through the whole paper, $\rm m_R$ means the 
maximal ideal of the local ring R. 

An element in
$\rm  Def(E_0, \alpha_0)(A')
\times_{ Def(E_0,\alpha_0)(A)} 
Def( E_0, \alpha_0 )(B)$ 
is given by a pair 
$\rm \left( \left( E',\alpha ' \right) , \left( F,\beta \right) \right)$ 
with
$\rm \left( E',\alpha ' \right) \in Def( E_0, \alpha_0)(A')$ 
and $\rm \left( F,\beta \right) \in Def( E_0, \alpha_0)(B)$, 
such that we have 
$\rm E'\otimes_{A'}A=F\otimes_{B}=E$ and 
$\rm \alpha '\otimes 1_A=\beta \otimes
1_A=\alpha$, where 
$\rm \left( E,\alpha \right) \in  Def( E_0, \alpha_0)(A)$. 
With $\rm A'\rightarrow A$ small also 
$\rm B'=B\times_{A}A'\rightarrow B$ is a small extension. 
Moreover there is an element $\rm t'\in B'$, such that
$\rm B=B'/t'B'$, $\rm t'\cdot m_{B'}=0$ and
$\rm B' \rightarrow A'$ induces $\rm t'B'\cong tA'$. 

By taking 
$\rm F'=E'\times_E F$ we obtain a commutative diagramm
\[
\begin{array}{ccccccccc}
\rm 0&\rightarrow&\rm E_0&
\stackrel{\rm t}{\longrightarrow}&
\rm E'&
\stackrel{\rm p}{\longrightarrow}&
\rm E & \rightarrow & 0 \\
&&\:\big\uparrow{\ssty=}&&\:\big\uparrow{\ssty\rm g}&&
\:\big\uparrow{\ssty\rm f}&&\\
\rm 0&\rightarrow&\rm E_0&
\stackrel{\rm t'}{\longrightarrow}&
\rm F'&
\stackrel{\rm q}{\longrightarrow}&
\rm F&\rightarrow&0
\end{array}
\]
with exact lines.
Since the induced mappings $\rm \bar{f}:\,F/m_BF\rightarrow E/m_AE$ and
$\rm \bar{f}:\,E'/m_{A'}E'\rightarrow E/m_AE$ are isomorphisms and
since $\rm t:\,E'/m_{A'}E'\rightarrow tE'$ with 
$\rm \left( x\; mod\; m_{A'}E' \right) \rightarrow tx$ is surjective,
we have
$\rm t'F'\cong tE'\cong tA'\otimes_{A'}E'\cong E_0\cong F/m_BF\cong t'B'\otimes_{B'}F'$, 
where the second isomorphism is given by the
$\rm A'$-flatness of $\rm E'$. 
Hence, with F is flat over B also 
$\rm F'$ is flat over $\rm A'$. 

By restriction we obtain
\[
\begin{array}{ccccccccc}
\rm 0&\rightarrow&\rm W_0&
\stackrel{\rm t}{-\!\!\!-\!\!\!\longrightarrow}&
\rm A'\otimes W_0&
-\!\!\!-\!\!\!-\!\!\!-\!\!\!-\!\!\!-\!\!\!-\!\!\!-\!\!\!\longrightarrow&
\rm A\otimes W_0 & \rightarrow & 0 \\
&&\:\Big\uparrow{\ssty =}&&
\hspace{5em}\Big\uparrow\rm\ssty\alpha'\circ g|_{SpecB'\times Y_0}&&
\Big\uparrow &&\\
0&
\rightarrow&
\rm W_0&
\stackrel{\rm t'\circ \alpha_0^{-1}}{-\!\!\!-\!\!\!\longrightarrow}&
\rm F'|_{SpecB'\times Y_0}&
\stackrel{\rm \beta\circ q |_{SpecB'\times Y_0}}
{-\!\!\!-\!\!\!-\!\!\!-\!\!\!-\!\!\!-\!\!\!-\!\!\!-\!\!\!\longrightarrow}&
\rm B\otimes W_0&
\rightarrow&
0
\end{array}
\]

For h the projection 
$\rm A'\otimes W_0 \rightarrow W_0$ the lower line splits by 
$\rm h\circ \alpha`\circ g|_{Y_0}$. 
We consider the composition
$\rm \beta '=\left( h\circ \alpha ' \circ g|_{SpecB'\times Y_0} \right) \oplus\left(  \beta \circ q|_{SpecB'\times Y_0} \right) : F'|_{SpecB'\times Y_0}\cong  B'\otimes W_0$
and obtain 
$\rm ( F',\beta )$ as an element in $\rm Def( E_0,\alpha _0 )( B')$.
The restriction of $\rm F'$ to $\rm SpecB\times X_0$ yields F, the pullback
of $\rm F'$ to $\rm SpecA'\times X_0$ yields $\rm E'$ and both is compatible with the framings. 
Thus the mapping 
$\rm ( ( E',\alpha ' ) ,( F,\beta )) \rightarrow ( F',\beta ' )$ 
is the inverse to 
\[ \rm Def( E_0, \alpha_0 )(B') \longrightarrow 
Def( E_0, \alpha_0)( A')\times_{ Def( E_0, \alpha_0)(A )} 
Def( E_0, \alpha_0)(B).\;\Box \]

\subsection{Deformations and extensions}

\begin{theorem}
\parbox[t]{30em}{
Let $\rm (E_0,\alpha_0)\in{\cal M}(X_0,Y_0,W_0)(\CC)$. 
There is a canonical isomorphism
$\rm Ext^1_{X_0}\left( E_0,E_0\otimes_{{\cal O}_{X_0}}{\cal J}_{Y_0|X_0}\right)
=Def \left( E_0,\alpha_0\right) \left( \CC [\varepsilon] \right)$.}
\end{theorem}

{\it Proof:}
For $\rm X=Spec\CC \left[ \varepsilon \right] \times X_0$ we have
$\rm {\cal O}_X= {\cal O}_{X_0}\oplus \varepsilon{\cal O}_{X_0}$. 
For
$\rm \left( E,\alpha, \psi\right)\in Def\left( E_0, \alpha{}_0 \right)
\left( \CC \left[ \varepsilon \right]  \right)$
the exact sequence
$\rm 0\rightarrow \varepsilon E \rightarrow E \rightarrow E/\varepsilon E \rightarrow 0$ yields therefore an element
\[ \rm \left( e \right)  = \left( 0\rightarrow E_0 
\stackrel{i}{\longrightarrow} E 
\stackrel{p}{\longrightarrow} E_0
\rightarrow 0
\right)\]
in $\rm Ext^1_{X_0}\left( E_0,E_0 \right)$,
where we have used $\rm \psi : E/\varepsilon E \cong E_0$ and
$\rm \CC \left[ \varepsilon \right] / \varepsilon \CC \left[ \varepsilon \right] \cong \varepsilon \CC \left[ \varepsilon \right] $.
By restriction to $\rm Y_0$ 
we obtain a commutative diagram with exact lines:
\[
\begin{array}{cccccccccc}
&
0& 
\rightarrow & 
\rm E_0|_{Y_0} &
\longrightarrow &
\rm E|_{Y_0} &
\longrightarrow &
\rm E_0|_{Y_0} &
\rightarrow & 
0\\
\rm \left( d \right) \hspace{1cm}&
&
&
\:{\big\downarrow} \ssty\alpha _0 &
&
{\big\downarrow} \ssty\alpha &
&
\:{\big\downarrow} \ssty\alpha _0 &
&
\\
&
0& 
\rightarrow & 
\rm W_0 &
\longrightarrow &
\rm W_0 \oplus \varepsilon W_0 &
\longrightarrow &
\rm W_0 &
\rightarrow 
& 0
\end{array}
\]
We note that
$\rm \left( e \right) \in ker \left( Ext^1_{X_0}\left( E_0,E_0 \right)
\rightarrow  Ext^1_{Y_0}\left( E_0|_{Y_0},E_0|_{Y_0} \right) \right) $.
Conversely, the pair ((e),(d)) determines the deformation
$\rm \left( E,\alpha, \psi\right) $ by defining the multiplication of 
a local section $\rm x \in E$ with $\varepsilon $ by
$\rm \varepsilon x = \left( i \circ p \right) \left( x \right) $.

Now we desribe a canonical map $\nu$ from
$\rm Ext^1_{X_0}\left( E_0,E_0\otimes_{{\cal O}_{X_0}}{\cal J}_{Y_0|X_0}\right)$ to
$\rm Def \left( E_0,\alpha_0\right) \left( \CC [\varepsilon] \right)$.
Let $\rm (f)=\left( 0\rightarrow E_0\otimes_{{\cal O}_{X_0}}
{\cal J}_{Y_0|X_0} \stackrel{i}{\longrightarrow} F \longrightarrow E_0 \rightarrow 0 \right)$. 
For j the embedding 
$\rm  E_0\otimes_{{\cal O}_{X_0}}{\cal J}_{Y_0|X_0}
\stackrel{j}{\hookrightarrow} E_0$ the fibre sum of i and j
\[
\begin{array}{ccccccccc}
0&
\rightarrow&
\rm E_0\otimes_{{\cal O}_{X_0}}{\cal J}_{Y_0|X_0}&
\rm \stackrel{i}{\longrightarrow}&
\rm F&
\longrightarrow &
\rm E_0 &
\rightarrow&
0 \\
&
&
\rm {\big \downarrow} \ssty j&
&
{\big \downarrow}&
&
{\big \downarrow}&
&
\\
0&
\rightarrow&
\rm E_0&
\rm {\longrightarrow}&
\rm E=F\amalg_{E_0\otimes {\cal J}_{Y_0|X_0}}E_0&
\longrightarrow &
\rm E_0 &
\rightarrow&
0
\end{array}
\]
defines an element 
$\rm (e)=(0\rightarrow E_0 \rightarrow E \rightarrow E_0 \rightarrow 0)$ in
$\rm Ext^1_{X_0}\left( E_0,E_0 \right)$.
In particular we have
$\rm E=F\amalg_{E_0\otimes {\cal J}_{Y_0|X_0}}E_0=\left( F\oplus E_0 \right) / N$,
where N is the subsheaf consisting of all pairs (x,y) such that there is a
$\rm z\in E_0\otimes {\cal J}_{Y_0|X_0}$ with i(z)=y and j(z)=x.
Thus we have 
\[ \rm E|_{Y_0}=E/{\cal J}_{Y_0|X_0}E= (F\oplus E_0)/
\left( N+( {\cal J}_{Y_0|X_0}F \oplus {\cal J}_{Y_0|X_0}E_0 ) \right)\]
and together with
$\rm N+{\cal J}_{Y_0|X_0}E_0 =i\left( {\cal J}_{Y_0|X_0}E_0 \right)
\oplus {\cal J}_{Y_0|X_0}E_0$
we obtain the composition
\[
\begin{array}{ll}
\rm E|_{Y_0}&
\rm =\left( F{\Big /}i\left( {\cal J}_{Y_0|X_0}E_0\right)+{\cal  J}_{Y_0|X_0}F \right)
\rm    \oplus \left( E_0 {\Big /} {\cal J}_{Y_0|X_0}E_0 \right)\\
&
\rm \cong \left( \left( F / i\left( {\cal J}_{Y_0|X_0}E_0\right) \right) 
      {\Big /}
\rm      {\cal J}_{Y_0|X_0}\left( F / i\left( {\cal J}_{Y_0|X_0}E_0 \right) \right) \right) \rm   \oplus \left( E_0 {\Big /} {\cal J}_{Y_0|X_0}E_0 \right)\\
&
\rm \cong \left( E_0 {\Big /} {\cal J}_{Y_0|X_0}E_0 \right) 
\rm     \oplus \left( E_0 {\Big /} {\cal J}_{Y_0|X_0}E_0 \right)\\
& 
\rm \cong E_0|_{Y_0}\oplus E_0|_{Y_0}
\begin{array}{c}
    \cong \\[-1.4em]
    -\!\!\!-\!\!\!-\!\!\!\longrightarrow\\[-1em]
    \ssty (\alpha_0, \varepsilon\alpha_0)
\end{array}
\rm W_0\oplus \varepsilon W_0,
\end{array}
\]
which yields a framing $\rm \alpha : E|_{Y_0} \cong W_0 \oplus \varepsilon W_0$
in a canonical way, whence we have found our map by taking
$\rm \nu\left( f \right)=\left( E,\alpha \right)$.

To find the inverse map of $\nu $, we associate for a given
$\rm \left( E,\alpha \right)\! \in \!
Def\! \left( E_0,\alpha_0\right) \left( \CC [\varepsilon] \right)$
a pair ((e),(d)) with $\rm (e)\in Ext^1_{X_0}\left( E_0,\alpha_0\right)$
as described above and denote with q the composition
\[\rm E \rightarrow E|_{Y_0} 
\stackrel{\alpha}{\longrightarrow} W_0 \oplus \varepsilon W_0
\stackrel{2nd\; projection}
{-\!\!\!-\!\!\!-\!\!\!-\!\!\!-\!\!\!-\!\!\!-\!\!\!-\!\!\!\longrightarrow} 
\varepsilon W_0 \cong W_0.\]
By defining 
$\rm F= ker(q)$
we obtain 
$\rm {\cal J}_{Y_0|X_0}E_0 \cong E_0 \times_E F$
and moreover, the map 
$\rm {\cal J}_{Y_0|X_0}E_0 \rightarrow F$
given by the fibre product is injective and its cokernel is isomorphic
to $\rm E_0$.
Thus we have obtained an element
$\rm (f)=\left( 0\rightarrow {\cal J}_{Y_0|X_0}E_0 \rightarrow F \rightarrow E_0 \rightarrow 0\right)$
in $\rm Ext^1_{X_0}\left( E_0,{\cal J}_{Y_0|X_0}E_0 \right)$
and a commutative diagram
\[
\begin{array}{ccccccccc}
&&0&&0&&&&\\
&&\downarrow&&\downarrow&&&&\\
\rm 0& \rightarrow &\rm {\cal J}_{Y_0|X_0}E_0 & \longrightarrow &\rm F& 
\rm \longrightarrow&\rm E_0& \rightarrow& 0\\
&&\big\downarrow&&\big\downarrow&&\hspace{1em}\big\downarrow \ssty =&&\\
\rm 0& \rightarrow &\rm E_0 & \longrightarrow & \rm E& 
\rm \longrightarrow&\rm E_0& \rightarrow& 0\\
\rm &&\big\downarrow&&\hspace{1ex}\big\downarrow\rm q&&&&\\
\rm &&\rm E_0|_{Y_0} & \stackrel{\alpha_0}{\longrightarrow}&\rm W_0&&&&\\
&&\downarrow&&\downarrow&&&&\\
&&0&&0&&&&
\end{array}
\]
with exact lines and columns,
from where we see that $\rm \nu (f) = (E,\alpha)$
and hence $\nu $ is a canonical bijection. 

It is not hard but some writing 
to describe the linear structure of 
$\rm Def ( E_0,\alpha_0)( \CC [\varepsilon])$
in terms of $\rm ((e),(d))$. 
From there it is obvious 
that $\nu$ is linear, too. $\Box$

\subsection{Smoothness of the deformation functor}
As in \cite{schlessinger} or \cite{artin2},
the functor
$\rm Def(E_0,\alpha_0)$
is called {\it formally smooth},
if for every small extension
$\rm R\rightarrow R/tR =\bar{R}$ of Artin rings (i.e. $\rm t\in R$ 
with $\rm t\cdot m_R=0$)
the induced mapping from
$\rm Def(E_0,\alpha_0)(R)$ to $\rm Def(E_0,\alpha_0)(\bar{R})$
is surjective.
In the following
we fix an arbitrary small extension
$\rm R\rightarrow R/tR =\bar{R}$ of Artin rings and define spaces
$\rm X=SpecR\times X_0$, $\rm Y=SpecR\times Y_0$, 
$\rm \bar{X}=Spec\bar{R}\times X_0$ and 
$\rm \bar{Y}=Spec\bar{R}\times Y_0$.
With W and $\rm \bar{W}$ we denote the pullbacks of $\rm W_0$
to Y and $\rm \bar{Y}$, respectively.
The assertion shown by the next Proposition is known as
{\it $\rm T^1$-lifting property}.

\begin{proposition}
Let $\rm ({E},{\alpha})\in Def(E_0,\alpha_0)({R})$ and
$\rm (\bar{E},\bar{\alpha})\in Def(E_0,\alpha_0)(\bar{R})$.
We may think of both as framed vector bundles to the framing data
$\rm (X,Y,W)$ and $\rm (\bar{X},\bar{Y},\bar{W})$, respectively.
If $\rm Ext^2_{X_0}(E_0,{\cal J}_{Y_0|X_0}\otimes E_0)=0$,
then the natural homomorphism between the deformation-theoretical
tangent spaces
$\rm Def(E,\alpha)(\CC[\varepsilon])\rightarrow
Def(\bar{E},\bar{\alpha})(\CC[\varepsilon])$ is surjective.
\end{proposition}

{\it Proof:}
With E flat over X and deformation of $\rm E_0$ we have
$\rm tE = tR\otimes E \cong E/m_RE = E_0$ and therefore a 
short exact sequence
\[ \rm 0 \rightarrow E_0 \stackrel{t}{\longrightarrow} E
\longrightarrow \bar{E} \rightarrow 0. \]
With $\rm \bar{E}$ flat over $\rm X_0$ and 
$\rm {\cal J}_{Y|X}={\cal J}_{Y_0|X_0}\otimes {\cal O}_X$ 
we also have
\[ \rm 0 \rightarrow {\cal J}_{Y|X}\otimes E_0 
\stackrel{t}{\longrightarrow} 
{\cal J}_{Y|X}\otimes E
\longrightarrow {\cal J}_{Y|X}\otimes \bar{E} \rightarrow 0. \]
As part of the long Ext-sequence we obtain
\[ \rm Ext^1_X(E,{\cal J}_{Y|X}\otimes E) \longrightarrow
Ext^1_{\bar{X}}(\bar{E},{\cal J}_{\bar{Y}|\bar{X}}\otimes \bar{E}) 
\longrightarrow
Ext^2_{X_0}(E_0,{\cal J}_{Y_0|X_0}\otimes E_0), \]
where 
$\rm Ext^1_X(E,{\cal J}_{Y|X}\otimes \bar{E})=
Ext^1_{\bar{X}}(\bar{E},{\cal J}_{\bar{Y}|\bar{X}}\otimes \bar{E})$
and 
$\rm Ext^2_X(E,{\cal J}_{Y|X}\otimes E_0)=
Ext^2_{X_0}(E_0,{\cal J}_{Y_0|X_0}\otimes E_0)$
were used.
The first homomorphism of the last exact sequence coincides with
$\rm Def(E,\alpha)(\CC[\varepsilon])\rightarrow
Def(\bar{E},\bar{\alpha})(\CC[\varepsilon])$
by Theorem 1.2 and we are done. $\Box$

\begin{theorem}
If $\rm dim_{\CC}Ext^1_{X_0} (E_0,{\cal J}_{Y_0|X_0}E_0)<\infty$ and
$\rm Ext^2_{X_0} (E_0,{\cal J}_{Y_0|X_0}E_0)=0$,
then $\rm Def(E_0,\alpha_0)$ is formally smooth.
\end{theorem}

{\it Proof:}
Because of Lemma 1.1 and since the dimension of
$\rm Ext^1_{X_0} (E_0,{\cal J}_{Y_0|X_0}E_0)$ is finite, 
the functor $\rm Def(E_0,\alpha_0)$ is prorepresentable by the criterion
of Schlessinger \cite{schlessinger}. 
The $\rm T^1$-lifting property holds by Proposition 1.3
and the application of the Kawamata-Ran principle
(\cite{kawamata}, Theorem 1) yields our assertion. $\Box$

{\it Remark.} D. Huybrechts and M. Lehn consider in \cite{huy-lehn1}
a framed module over an algebraic nonsingular
projective variety $\rm X_0$ given by a coherent 
$\rm {\cal O}_X$-module $\rm E_0$ and a nonzero morphism 
from $\rm E_0$ to a fixed coherent
$\rm {\cal O}_X$-module $\rm D_0$. 
They introduce a notion
of stability in this situation and describe the 
deformation-theoretical tangent space of a stable framed module
$\rm (E_0, E_0 \rightarrow D_0)$ as 
$\rm \EE xt^1(E_0, E_0 \rightarrow D_0)$
and the obstructions for smoothness as living in
$\rm \EE xt^2(E_0, E_0 \rightarrow D_0)$.
Here we think of a sheaf as a complex of sheaves concentrated in 
zero and of a morphism of sheaves as a complex of sheaves concentrated
in zero and one.
From the short exact sequence of complexes
\[ \rm 0 \rightarrow ({\cal J}_{Y_0|X_0}E_0 \rightarrow 0)
\rightarrow (E_0 \rightarrow E_0|_{Y_0}) \rightarrow
(E_0|_{Y_0} \rightarrow E_0|_{Y_0}) \rightarrow 0 \]
we obtain the long hyperext-sequence
\[ 
\begin{array}{lr}
\rm 
\ldots \rightarrow 
\EE xt^{i-1}(E_0, E_0|_{Y_0} \rightarrow E_0|_{Y_0})
\rightarrow 
\EE xt^i(E_0, {\cal J}_{Y_0|X_0}E_0))
\rightarrow & \\ & \rm
\hspace{-14em}
\EE xt^i(E_0,E_0 \rightarrow E_0|_{Y_0}) 
\rightarrow
\EE xt^i(E_0, E_0|_{Y_0} \rightarrow E_0|_{Y_0})
\rightarrow \ldots 
\end{array}
\]
Since 
$\rm Hom^{\bullet}(E_0,(E_0|_{Y_0} \rightarrow E_0|_{Y_0}))$
is the exact complex 
$\rm Hom(E_0,E_0|_{Y_0})\rightarrow Hom(E_0,E_0|_{Y_0})$
we have
$\rm \EE xt^i(E_0, E_0|_{Y_0} \rightarrow E_0|_{Y_0})
=H^i(\RR^{+}Hom^{\bullet}(E_0,(E_0|_{Y_0} \rightarrow E_0|_{Y_0})))=0$ 
(cf. \cite{hartshorne2},\S 1.6).
With 
$\rm \EE xt^i(E_0, {\cal J}_{Y_0|X_0}E_0))=
Ext^i(E_0, {\cal J}_{Y_0|X_0}E_0))$
the results above correspond to the results in \cite{huy-lehn1}
in the smooth projective and stable case.

\section{Framed instanton bundles}

\subsection{The restriction 
        $\rm {\cal M}(P,F,{\cal O}_F^r)\rightarrow {\cal M}(S,F,{\cal O}_F^r)$}
Let $\rm X_0$, $\rm Y_0$, $\rm W_0$ and $\rm {\cal M}(X_0,Y_0,W_0)$ be as
in section 1.1. Then the framing data $\rm (X_0,Y_0,W_0)$ are called
{\it simplifying}, if for any two representatives
$\rm (E_1,\alpha_1)$, $\rm (E_2,\alpha_2)$ of elements in 
$\rm {\cal M}(X_0,Y_0,W_0)(\CC)$ we have
$\rm H^0(X_0,\Hom(E_1,E_2)\otimes{\cal J}_{Y_0|X_0})=0$.
This definition of simplifying is analogous to the definition in
\cite{lehn} or \cite{lubke}.

We consider now a twistor fibration
$\rm \pi : P\rightarrow M$ with M compact and
$\rm S\subset P$ a surface of degree 1 on the 
twistor fibres containing a unique twistor fibre 
$\rm F=S\cdot \bar{S}\subset S$ as in the introduction.
Moreover we fix some integer $\rm r>0$.
\begin{lemma}
{The framing data $\rm (S,F,{\cal O}_F^r)$ are simplifying.}
\end{lemma}
{\it Proof:} Consider 
$\rm (E_1,\alpha_1),(E_2,\alpha_2)\in {\cal M}(S,F,{\cal O}_F^r)(\CC)$
and let s be in
$\rm H^0(S,\Hom(E_1,E_2)(-F))$. 
If H denotes the Douady space of lines on S, then there is an open neighbourhood
$\rm F\in U$ in H, such that $\rm E_1$ and $\rm E_2$ and hence
$\rm \Hom(E_1,E_2)$ are trivial along all lines $\rm L\in U$
and moreover $\rm (F\cdot L)=1$. 
Since s is a global section of $\rm \Hom(E_1,E_2)$ vanishing
along F, for all $\rm L\in U$ we have $\rm s|_L$ as global section of
a trivial vector bundle over a line vanishing in one point and therefore
$\rm s|_L=0$. Thus s vanishes on an open subset of S and with S irreducible
we obtain, that s vanishes everywhere, which proves our claim. $\Box$
\begin{lemma}
If $\rm W$ is a vector bundle on $\rm S$ with $\rm W|_F\cong {\cal O}_F^r$,
then the framing data $\rm (P,S,W)$ are simplifying.
\end{lemma}
{\it Proof:} Consider $\rm (E_1,\alpha_1),(E_2,\alpha_2)\in {\cal M}(P,S,W)(\CC)$,
$\rm s\in H^0(P,\Hom(E_1,E_2)(-S))$ and let now H denotes 
the Douady space of lines on P. We find an open neighbourhood 
$\rm F\in U$ in H, such that 
$\rm \Hom(E_1,E_2)$ is trivial along all lines $\rm L\in U$
and moreover $\rm (S\cdot L)=1$. Since s vanishes along S, it vanishes
also along all L in U, therefore on an open subset of P and therefore
everywhere. $\Box$

Finally we also state
\begin{lemma}
The framing data $\rm (P,F,{\cal O}_F^r)$ are simplifying.
\end{lemma}
{\it Proof: } Let $\rm (E_1,\alpha_1),(E_2,\alpha_2)\in {\cal M}(S,F,{\cal O}_F^r)(\CC)$
and $\rm s\in H^0(P,\Hom(E_1,E_2)\otimes{\cal J}_{F|P})$.
Then by lemma 2.1 s vanishes along S and thus everywhere on P by
lemma 2.2. $\Box$\vspace{1ex}

\noindent
{\it Remark:}
If the framing data $\rm (X_0,Y_0,W_0)$ are simplifying, then 
we have for all elements
$\rm (E,\alpha)$
in $\rm {\cal M}(X_0,Y_0,W_0)(\CC)$ in particular
$\rm Aut((E,\alpha))=\{ id_E\}$. Because, for 
$\rm s\in Aut((E,\alpha))$ we have $\rm s|_{Y_0}=id_E|_{Y_0}$ and therefore
$\rm s-id_E\in H^0(X_0,\Hom(E,E)\otimes{\cal J}_{Y_0|X_0})=0$.
Since by \cite{poon} F intersects every effective divisor on P positively,
every morphism of framed vector bundles to the framing data
$\rm (P,F,{\cal O}_F^r)$ is an isomorphism and moreover there is at most
one morphism between two such framed vector bundles.
\indent

\begin{proposition}
For $\rm {\cal M}(P,F,{\cal O}_F^r)\rightarrow {\cal M}(S,F,{\cal O}_F^r)$ the 
natural transformation given by restriction, the fibre over
a point 
$\rm (W,\sigma)\in {\cal M}(S,F,{\cal O}_F^r)(\CC )$ is naturally equivalent
to $\rm {\cal M}(P,S,W)$.
\end{proposition}
{\it Proof:} 
Let X be some analytic space and $\rm V_1$ and $\rm V_2$ 
two vector bundles on the space $\rm X\times P$ with
$\rm V_1|_{X\times F}\cong V_2|_{X\times F}\cong {\cal O}_{X\times F}^r$.
The exactness of
\[ \rm H^0(S,\Hom(V_1|_S,V_2|_S)(-X\times F)) \rightarrow
       Hom(V_1|_{X\times S},V_2|_{X\times S}) 
      \rightarrow Hom(V_1|_{X\times F},V_2|_{X\times F})  \]
and lemma 2.1 yield the injectivity of 
\[ \rm Hom(V_1|_{X\times S},V_2|_{X\times S}) 
       \rightarrow Hom(V_1|_{X\times F},V_2|_{X\times F}),  \]
a fact which will be denoted with $\rm (\ast)$ in the following.

For W  a vector bundle on S and $\rm \sigma : W|_F \cong {\cal O}_F^r$
we consider the natural transformation 
$\rm {\cal M}(P,S,W)\rightarrow {\cal M}(P,F,{\cal O}_F^r)$ given by restriction
of the framing map.
I.e. if $\rm (V,\alpha)$ is an element in
$\rm {\cal M}(P,S,W)(X)$, then $\rm (V,\alpha)$ is associated to the element
$\rm (V,(p^{\ast}_F\sigma) \circ (\alpha |_{X\times F}))$, where the
morphism
$\rm p_F: X\times F \rightarrow F$ is the projection and
$\rm p^{\ast}_F\sigma: p^{\ast}_FW\cong {\cal O}^r_{X\times F}$ is the pullback.
If $\rm (V_1,\alpha_1)$ and $\rm (V_2,\alpha_2)$ are in
$\rm {\cal M}(P,S,W)(X)$, such that there is an isomorphism of framed vector bundles
\[ \rm \varphi : (V_1,(p^{\ast}_F\sigma) \circ (\alpha_1|_{X\times F}))
\cong (V_2,(p^{\ast}_F\sigma) \circ (\alpha_2|_{X\times F})), \]
i.e. with
$\rm (p^{\ast}_F\sigma) \circ (\alpha_2|_{X\times F})\circ (\varphi |_{X\times F}) = (p^{\ast}_F\sigma) \circ (\alpha_1|_{X\times F})$,
then $\rm (\ast)$ implies, that 
$\rm \varphi |_{X\times S}= \alpha_2^{-1}\circ \alpha_1$ and therefore
$\rm (V_1,\alpha_1)\cong (V_2,\alpha_2)$. Thus 
$\rm {\cal M}(P,S,W)(X) \rightarrow {\cal M}(P,F,{\cal O}_F^r)(X)$ is injective.

Now we consider an element 
$\rm (V,\alpha )\in {\cal M}(P,F,{\cal O}_F^r)(X)$, such that there is an isomorphism
of framed vector bundles
$\rm \varphi : (V|_{X\times S},\alpha )\cong (p^{\ast}_SW, p^{\ast}_F\sigma )$
with $\rm p_S: X\times S \rightarrow S$ the projection.
With $\rm (\ast)$ and with 
$\rm (p^{\ast}_F\sigma)\circ (\varphi |_{X\times F})=\alpha$,
$\rm \varphi$ is well defined by 
$\rm (p^{\ast}_F\sigma)^{-1}\circ \alpha \in 
 Hom(V|_{X\times F},p^{\ast}_FW)$
and hence $\rm (V,\alpha )$ uniquely determines an element 
$\rm (V,\varphi ) \in {\cal M}(P,S,W)(X)$.
Conversely, the image of $\rm (V,\varphi )$ in 
$\rm {\cal M}(P,F,{\cal O}_F^r)(X)$ is again $\rm (V,\alpha )$,
hence the natural transformation
$\rm {\cal M}(P,S,W)\rightarrow {\cal M}(P,F,{\cal O}_F^r)$ induces a bijection between
$\rm {\cal M}(P,S,W)(X)$ and the fibre of
$\rm {\cal M}(P,F,{\cal O}_F^r)(X)\rightarrow {\cal M}(S,F,{\cal O}_F^r)(X)$
over $\rm (p^{\ast}_SW,p^{\ast}_F\sigma)$.

These bijections correspond actually to a natural equivalence between
$\rm {\cal M}(P,S,W)$  and the fibre of
$\rm {\cal M}(P,F,{\cal O}_F^r)\rightarrow {\cal M}(S,F,{\cal O}_F^r)$
over $\rm (W,\sigma)$. So far everything was natural except the choice
of $\rm (W,\sigma)$ as a representative of an element in
$\rm {\cal M}(S,F,{\cal O}_F^r)(\CC )$. If $\rm (W',\sigma')$ is another
representative of the same isomorphism class, then there is an
isomorphism 
$\rm \psi : W' \rightarrow W$ with
$\rm \sigma \circ \psi |_F = \sigma$. 
Because of $\rm (\ast)$ this isomorphism is unique and therefore there is
a unique natural equivalence from
$\rm {\cal M}(P,S,W')$ to $\rm {\cal M}(P,S,W)$, which moreover maps the fibre of
$\rm {\cal M}(P,F,{\cal O}_F^r)\rightarrow {\cal M}(S,F,{\cal O}_F^r)$ over $\rm (W,\sigma)$ 
identically to itself. $\Box$

Since Chern classes are locally constant, the functors of families
of framed vector bundles 
$\rm {\cal M}(P,F,{\cal O}_F^r)$, $\rm \:{\cal M}(S,F,{\cal O}_F^r)$ and 
$\rm \,{\cal M}(P,S,W)$ as above split into open and closed subfunctors
$\rm {\cal M}(P,F,{\cal O}_F^r,c_{\bullet})$, 
$\rm {\cal M}(S,F,{\cal O}_F^r,c_{\bullet})$ and
$\rm {\cal M}(P,S,W,c_{\bullet})$, 
where $\rm c_{\bullet}$ are fixed
Chern classes in $\rm H^{\bullet}(P,\ZZ)$ and $\rm H^{\bullet}(S,\ZZ)$,
respectively.
In particular, our restriction 
$\rm {\cal M}(P,F,{\cal O}_F^r)\rightarrow {\cal M}(S,F,{\cal O}_F^r)$
splits into open and closed parts
\[ \rm {\cal M}(P,F,{\cal O}_F^r,c_{\bullet})\rightarrow 
{\cal M}(S,F,{\cal O}_F^r,i^{\ast}c_{\bullet})\] with
$\rm c_{\bullet}\in H^{\bullet}(P,\ZZ)$ and 
$\rm i: S \hookrightarrow P$.
It is easy to see, that the fibre of
$\rm {\cal M}(P,F,{\cal O}_F^r)\rightarrow {\cal M}(S,F,{\cal O}_F^r)$
over an 
$\rm (W,\sigma)\in {\cal M}(S,F,{\cal O}_F^r,i^{\ast}c_{\bullet})(\CC)$
is just $\rm {\cal M}(P,S,W,c_{\bullet})$.
Later on, when we restrict ourself to the case of framed
instanton bundles, we only have to consider Cern classes coming from
$\rm H^{\bullet}(M,\ZZ)$, i.e. we will consider
\[ \rm {\cal M}(P,F,{\cal O}_F^r,\pi^{\ast}c_{\bullet})\rightarrow 
{\cal M}(S,F,{\cal O}_F^r,\pi^{\ast}_Sc_{\bullet}), \]
where $\rm c_{\bullet}\in H^{\bullet}(M,\ZZ)$ and 
$\rm \pi_S$ is the restriction of the twistor fibration to S.

The image of our natural transformation has a fine moduli space
by the following proposition:

\begin{proposition}
The functor $\rm {\cal M}(S,F,{\cal O}_F^r,c_{\bullet})$ is represented
by a separated algebraic space of finite type over $\CC$.
\end{proposition}

{\it Proof:} This is a corollary of Satz 3.4.1 in Lehn's \cite{lehn}.
By \cite{kurke1}, Prop. 2.1, S is a smooth projective
surface and by Lemma 2.1 the framing data $\rm (S,F,{\cal O}_F^r)$
are simplifying. Therefore it is sufficient to show, that F is the
support of a divisor, which is big and nef.
Since the linear system $\rm |F|$ defines by \cite{kurke1} a blowing up
$\rm S\rightarrow \PP^2_{\CC}$, F intersects all lines on S coming from
$\PP^2_{\CC}$ with multiplicity 1 and all exceptional lines with multiplicity
0, which implies that F is numerically effective.
The big-condition follows by $\rm {\cal N}_{F|S}\cong {\cal O}_F(1)$. $\Box$

\subsection{Smoothness and dimensions}
By H. Kurke's article \cite{kurke1}, Prop.2.1, 
the canonical bundle of P is 
$\rm K_P={\cal O}_P(-2S-2\bar{S})$. 
Moreover, 
$\rm H^{\bullet}(P,\ZZ)$ is generated by elements
$\rm \{ \omega,\eta_1,\ldots ,\eta_n\}\in H^2(P,\ZZ)$,
where $\rm \{ \eta_1,\ldots ,\eta_n\}$
are an orthogonal basis of 
$\rm H^2(M,\ZZ)$.
We put $\rm \eta = \sum \eta_i$.
There are the relations 
$\rm \eta_i^2=-[F]$, $\rm \eta_i^3=0$, $\rm \eta_i\eta_j=0$
for $\rm i\not= j$, $\rm \omega^2 + \omega\eta = [F]$,
$\rm \omega^2 \eta_i =1$, $\rm \omega \eta_i^2 = -1$,
$\rm c_1(P)=4\omega + 2\eta$ and $\rm c_2(P)=3(e(M)-sgn(M))[F]$.
We have 
$\rm [S]=\omega + \sum a_i\eta_i$ with 
$\rm a_i=0\; or\; 1$ and from
$\rm (S\cdot \bar{S})=F$ we easily infer
$\rm [\bar{S}]=\omega + \sum (1-a_i)\eta_i$.
We put $\rm \sigma = \sum a_i\eta_i$ and 
$\rm \bar{\sigma}= \sum (1-a_i)\eta_i $.

We fix a Chern class 
$\rm c_{\bullet}\in H^{\bullet}(M,\ZZ)$ and consider a vector bundle V on P
of rank r
with $\rm c_{\bullet}(V)=\pi^{\ast}c_{\bullet}$.

\begin{lemma}
We have the equality of Euler characteristics:
$\rm \chi(\End (V)(-S))= \chi (\End (V)(-\bar{S}))$.
\end{lemma}

{\it Proof:}
By \cite{hirzebruch}, Theorem 4.4.3 we have
\[ \rm c_1(\End V)=c_3(\End V)=0 \; and\; 
c_2(\End V)=2rc_2(V)+(1-r)c_1^2(V).\]
Since the Chern classes of V are coming from $\rm H^{\bullet}(M,\ZZ)$,
we have $\rm c_2(\End V)$ as a multiple of $\rm [F]$.
From the same theorem we obtain by computation that
\[ \begin{array}{l}
\rm c_1(\End (V)(-S))=r(\omega + \sigma),\\
\rm c_2(\End (V)(-S))=c_2(\End V)+
       \displaystyle{r \choose 2}(\omega^2 + 2\omega\sigma +\sigma^2),\\
\rm c_3(\End (V)(-S))=(r-2)c_2(\End V)(\omega + \sigma )+
       \displaystyle{r \choose 3}
       (\omega^3 + 3\omega^2\sigma + 3\omega\sigma^2),\\
\rm c_1(\End (V)(-\bar{S}))=r(\omega + \bar{\sigma}),\\
\rm c_2(\End (V)(-\bar{S}))=c_2(\End V)+
       \displaystyle{r \choose 2}
       (\omega^2 + 2\omega\bar{\sigma} +\bar{\sigma}^2),\\
\rm c_3(\End (V)(-\bar{S}))=(r-2)c_2(\End V)(\omega + \bar{\sigma} )+
       \displaystyle{r \choose 3}
       (\omega^3 + 3\omega^2\bar{\sigma} + 3\omega\bar{\sigma}^2),\\.
\end{array} \]

The Hirzebruch-Riemann-Roch theorem yields
\[ \begin{array}{ll}
\rm \chi(\End (V)(-S))= &
\rm \displaystyle 
    \frac{r}{24}c_1(P)c_2(P)
     +\frac{r\omega + r\sigma}{12}(16\omega^2+4\eta^2 +16 \omega\eta 
     +3(e(M)-sgn(M))[F])
     \\[2ex] &
\rm \displaystyle 
      +\frac{2\omega + \eta}{2}\left( r^2\omega^2  + r^2\sigma^2 +2r^2\omega\sigma 
     -2c_2(\End V)- (r^2-r)(\omega^2 + \sigma^2 +2\omega\sigma)\right)
     \\[2ex] &
\rm \displaystyle 
     +\frac{r^3}{6}(\omega^3 + 3\omega^2\sigma + 3\omega\sigma^2)
     -\frac{r\omega+r\sigma}{2}(c_2(\End V) 
     + \frac{r^2-r}{2}(\omega^2 + \sigma^2 + 2\omega \sigma))
      \\[2ex] &
\rm \displaystyle 
     +\frac{r-2}{6}c_2(\End V)(\omega + \sigma)
     +\frac{r^3-3r^2+2r}{36}(\omega^3 + 3\omega^2 \sigma + 3\omega\sigma^2)
\end{array} \]
and the analogous formula for $\rm \chi(\End (V)(-\bar{S}))$, where
just all $\sigma$'s are replaced by $\bar{\sigma}$.
Thus we obtain 
\[ \begin{array}{lr}
\rm \chi(\End (V)(-S))- \chi (\End (V)(-\bar{S})) = &\\
&\rm\hspace{-12em}
      \omega(\sigma^2 - \bar{\sigma}^2)
      (\frac{7}{6}r+\frac{1}{2}r^2-\frac{1}{6}r^3)
     +\omega^2(\sigma - \bar{\sigma})
      (\frac{13}{6}r+\frac{1}{2}r^2-\frac{1}{6}r^3)
     +\omega\eta(\sigma - \bar{\sigma})r ,
\end{array} \]
where we have used 
\[ \rm (\sigma - \bar{\sigma})\eta^2=(\sigma^2 - \bar{\sigma}^2)\eta
=(\sigma - \bar{\sigma})[F]=(\sigma - \bar{\sigma})c_2(\End V)
=(\sigma - \bar{\sigma})(\omega^2+\omega\eta)=0 .\]
We had 
$\rm \sigma = \sum a_i\eta_i$ 
and $\rm \bar{\sigma} = \sum (1-a_i)\eta_i$
with $\rm a_i=0,1$ and put 
$\rm A=\sum a_i$.
Thus 
$\rm \omega(\sigma^2 - \bar{\sigma}^2)=n-2A$
and 
$\rm \omega^2(\sigma - \bar{\sigma})=2A-n$.
With $\rm \eta\sigma=\sigma^2$ and 
$\rm \eta\bar{\sigma}=\bar{\sigma}^2$ we also have
$\rm \omega\eta(\sigma - \bar{\sigma})=\omega(\sigma^2 - \bar{\sigma}^2)
=n-2A$.
Inserting all this in the equation above we obtain
\[ \rm \chi(\End (V)(-S))- \chi (\End (V)(-\bar{S})) = 0,\; \Box\]

Due to Hitchin's article \cite{hitchin} we have for all 
$\rm U(r)-$instanton bundles V on P the vanishing
$\rm H^1(P,V(-S-\bar{S}))=0$. If V is an $\rm U(r)-$instanton bundle,
then so does $\rm \End V$ and we have also
$\rm H^1(P,\End (V)(-S-\bar{S}))=0$.
We denote with $\rm {\cal M}_0(P,F,{\cal O}_F^r)$ the open subfunctor of
$\rm {\cal M}(P,F,{\cal O}_F^r)$ with 
\[ \rm {\cal M}_0(P,F,{\cal O}_F^r)(\CC)=\{ (V,\alpha)\in 
{\cal M}(P,F,{\cal O}_F^r)(\CC) | \, 
H^1(P,\End (V)(-S-\bar{S}))=0 \}. \]
Analogously to our notations in Section 2.1 we put
\[ \rm {\cal M}_0(P,F,{\cal O}_F^r,c_{\bullet})
 ={\cal M}_0(P,F,{\cal O}_F^r)\cap {\cal M}(P,F,{\cal O}_F^r,c_{\bullet}). \]

\begin{lemma}
\begin{tabular}[t]{l}
Consider $\rm c_{\bullet}\in H^{\bullet}(M,\ZZ)$ and
$\rm  (V,\alpha)\in 
      {\cal M}_0(P,F,{\cal O}_F^r,\pi^{\ast}c_{\bullet})(\CC)$.\\[1ex] 
Then we have
$ \rm dim_{\CC}H^i(P,\End (V)(-S))=\left\{
\begin{array}{ll}
\rm -\chi(\End(V)(-S))& \rm for\; i=1\\
\rm 0 & \rm else.
\end{array} \right. $
\end{tabular}
\end{lemma}

{\it Proof:}
By Lemma 2.3 we have 
\[ \rm H^0(P,\End (V)(-S))=0.\] 
With Serre duality and with $\rm K_P={\cal O}_P(-2S-2\bar{S})$ we obtain
\[ \rm H^3(P,\End (V)(-S))^{\vee}=
       H^0(P,\End (V)(-S-2\bar{S}))\subset
       H^0(P,\End (V)(-S)) =0 . \]

For $\rm \bar{E}$ be the exceptional divisor of the blowing up
$\rm |F|: \bar{S}\rightarrow \PP^2$ we have
$\rm K_{\bar{S}}={\cal O}_{\bar{S}}(-3F+\bar{E})$. With
$\rm K_{\bar{S}}={\cal O}_{\bar{S}}\otimes {\cal O}_P(\bar{S})\otimes K_P$
we infer the short exact sequence
\[ \rm 0 \rightarrow \End (V)(-S-2\bar{S}) \rightarrow 
      \End (V)(-S-\bar{S}) \rightarrow  \End (V|_{\bar{S}})(-2F+\bar{E})
       \rightarrow 0. \]
Since the bundle
$\rm \End (V|_{\bar{S}})$ is trivial on general lines of $\rm \bar{S}$
and $\rm \bar{E}$ is exceptional, also $\rm \End (V|_{\bar{S}})(\bar{E})$
is trivial on general lines and analogously to the proof of Lemma 2.1
we obtain the vanishing of
$\rm H^0(\bar{S},\End (V|_{\bar{S}})(-2F+\bar{E}))$.
Hence, the cohomology sequence of the last short exact sequence
implies 
an inclusion 
$\rm H^1(P, \End (V)(-S-2\bar{S}))\hookrightarrow 
        H^1(P,\End (V)(-S-\bar{S}))$
and our presumption gives us the vanishing
$\rm H^1(P, \End (V)(-S-2\bar{S}))=0 $.
Then Serre duality yields
\[ \rm H^2(P, \End (V)(-S))^{\vee}=H^1(P, \End (V)(-S-2\bar{S}))=0.\; \Box \]

By Theorem 1.4 the last lemma gives us:

\begin{corollary}
Around every closed point
$\rm (V,\alpha)\in {\cal M}_0(P,F,{\cal O}_F^r,\pi^{\ast}c_{\bullet})(\CC)$,
the fibres of the natural transformation
\[ \rm {\cal M}_0(P,F,{\cal O}_F^r,\pi^{\ast}c_{\bullet})
   \rightarrow {\cal M}(S,F,{\cal O}_F^r,\pi_S^{\ast}c_{\bullet}) \]
are smooth and of constant dimension
$\rm -\chi(\End(V)(-S))$.
\end{corollary}

\begin{theorem}
The natural transformation
\[ \rm {\cal M}_0(P,F,{\cal O}_F^r,\pi^{\ast}c_{\bullet})
   \rightarrow {\cal M}(S,F,{\cal O}_F^r,\pi_S^{\ast}c_{\bullet}) \]
is smooth.
More detailed we have the tangent space of the functor
$\rm {\cal M}_0(P,F,{\cal O}_F^r,\pi^{\ast}c_{\bullet})$ in a point
$\rm (V,\alpha)\in {\cal M}_0(P,F,{\cal O}_F^r,\pi^{\ast}c_{\bullet})(\CC)$
as a natural direct sum of the tangent space in direction of the fibre
and the tangent space in direction of the image, i.e.
\[ \rm T_{(V,\alpha)}{\cal M}_0(P,F,{\cal O}_F^r,\pi^{\ast}c_{\bullet})=
T_{(V,id_{V|_S})}{\cal M}(P,S,V|_S,\pi^{\ast}c_{\bullet})
\oplus
T_{(V|_S,\alpha)}{\cal M}(S,F,{\cal O}_F^r,\pi_S^{\ast}c_{\bullet}) .\]
Both summands are of dimension $\rm -\chi(\End(V)(-S))$. 

$\rm {\cal M}_0(P,F,{\cal O}_F^r,\pi^{\ast}c_{\bullet})$ and 
the fibres over 
$\rm {\cal M}(S,F,{\cal O}_F^r,\pi_S^{\ast}c_{\bullet})$
are smooth and equidimensional.
The image of the natural transformation in
$\rm {\cal M}(S,F,{\cal O}_F^r,\pi_S^{\ast}c_{\bullet})$
is an open and smooth algebraic subspace.
\end{theorem}

{\it Proof:}
The diagram $\rm (\ast)$
\[ \begin{array}{ccccccccc}
&&&&0&&0&&\\
&&&&\big\downarrow&&\big\downarrow&&\\
&&&&\rm {\cal O}_P(-S)&\rm =&\rm {\cal O}_P(-S)&&\\
&&&&\big\downarrow&&\big\downarrow&&\\
0&\longrightarrow&\rm {\cal O}_P(-S-\bar{S})&\longrightarrow&
 \rm {\cal O}_P(-S)\oplus{\cal O}_P(-\bar{S})&\longrightarrow&
 \rm {\cal J}_{F|P}&\longrightarrow&0\\
&&\big\downarrow =&&\big\downarrow&&\big\downarrow&&\\
0&\longrightarrow&\rm {\cal O}_P(-S-\bar{S})&\longrightarrow&
 \rm {\cal O}_P(-\bar{S})&\longrightarrow&
 \rm {\cal O}_S(-F)&\longrightarrow&0\\
&&&&\big\downarrow&&\big\downarrow&&\\
&&&&0&&0&&
\end{array} \]
is commutative with exact lines and columns.
For $\rm (V,\alpha)\in {\cal M}_0(P,F,{\cal O}_F^r,\pi^{\ast}c_{\bullet})(\CC)$
and because of the vanishing of the neighboured cohomology groups
we obtain
\[ \begin{array}{ccc}
0&&0\\
\downarrow & & \downarrow \\
\rm H^1(P,\End (V)(-S))&\stackrel{=}{\longrightarrow}&\rm H^1(P, \End (V)(-S))\\
\downarrow & & \downarrow \\
\begin{array}{c}
 \rm H^1(P,\End (V)(-S))\\[-1ex]
 \oplus\\[-1ex]
 \rm H^1(P,\End (V)(-\bar{S}))
\end{array}&
\stackrel{\cong}{\longrightarrow}&\rm H^1(P, \End (V)\otimes {\cal J}_{F|P})\\
\downarrow & & \downarrow \\
\rm H^1(P,\End (V)(-\bar{S}))&\stackrel{\cong}{\longrightarrow}&
 \rm H^1(S,\End (V|_S)(-F))\\
\downarrow & & \downarrow \\
0&&0
\end{array} \]
as a commutative diagram with exact columns.
The second column contains the tangent maps at $\rm (V,\alpha)$ 
corresponding to the composition
\[ \rm {\cal M}(P,S,V|_S,\pi^{\ast}c_{\bullet}) \hookrightarrow
  {\cal M}_0(P,F,{\cal O}_F^r,\pi^{\ast}c_{\bullet}) \rightarrow
  {\cal M}(S,F,{\cal O}_F^r,\pi_S^{\ast}c_{\bullet}) \]
Therefore the natural transformation 
\[ \rm 
  {\cal M}_0(P,F,{\cal O}_F^r,\pi^{\ast}c_{\bullet}) \rightarrow
  {\cal M}(S,F,{\cal O}_F^r,\pi_S^{\ast}c_{\bullet}) \]
is submersive and the tangent space of 
$\rm {\cal M}_0(P,F,{\cal O}_F^r,\pi^{\ast}c_{\bullet})$
in $\rm (V,\alpha)$
is a natural direct sum of the tangent spaces in direction of the fibre and in direction of the image.

With Lemma 2.6 we have
\[ \rm h^1\!(\End (V)(-S))=-\chi (\End (V)(-S)) =-\chi (\End (V)(-\bar{S})) 
   =h^1\!(\End (V)(-\bar{S}))=h^1(\End (V|_S)(-F)), \]
whence both direct summands are of the same dimension
$\rm -\chi (\End (V)(-S))$.
From $\rm (\ast)$ we also obtain 
\[ \rm H^2(P, \End (V)\otimes {\cal J}_{F|P})\cong 
  H^2(S,\End (V|_S)(-F)) \] and
\[ \rm H^2(S,\End (V|_S)(-F))\cong 
  H^3(P, \End (V)(-S-\bar{S})). \]
We have
\[ \rm H^3(P, \End (V)(-S-\bar{S}))^{\vee} =
       H^0(P, \End (V)(-S-\bar{S})) \subset
       H^0(P, \End (V)(-S-\bar{S}))=0 \]
and hence
$\rm H^2(S,\End (V|_S)(-F))=0 $ and 
$\rm H^2(P, \End (V)\otimes {\cal J}_{F|P})=0$,
which shows the smoothness by Theorem 1.4. $\Box$

\subsection{ Framed U(r)-instanton bundles}

\begin{lemma}
The property of a vector bundle on $\rm P$ to be a mathematical
instanton bundle is open.
\end{lemma}

{\it Proof:}
Let T be an analytical space and V a vector bundle over $\rm T\times P$.
For $\rm t\in T$ we denote with $\rm V_t$ the vector bundle on P
induced by $\rm V|_{\{ t \} \times P}$. We assume that for a point $\rm 0\in T$
the vector bundle $\rm V_0$ is an instanton bundle, 
i.e. trivial along twistor fibres.
We have to show that there is a neighbourhood $\rm 0\in T_1 \subset T$, such
that for all $\rm t\in T_1$ the vector bundle $\rm V_t$ is an instanton
bundle.

We consider the universal family 
\[ \begin{array}{ccc}
\rm Z & \stackrel{\nu}{\longrightarrow} & \rm P \\
\rm \big\downarrow \ssty\mu &&\\
\rm H &&
\end{array} \]
of lines on P as in the introduction. We may restrict H to the open neighbourhood
of M consisting of points, where the corresponding lines have intersection
product one with S. 
We put $\rm \bar{\nu}=id_T \times \nu$,
\linebreak[4] 
$\rm \bar{\mu}=id_T \times\mu$
and $\rm W=\bar{\nu}^{\ast} (V\otimes {\cal O}_{T\times P}(-\, T\times S))$.
Recall that $\mu$ is a proper and flat morphism with projective lines as fibres.
There are an covering
$\rm \{ U^{\lambda}_0|\, \lambda \in \; some\; index\; set \}$ of H of sufficiently  small open neighbourhoods and  effective divisors 
$\rm D^{\lambda}\subset\bar{\mu}^{-1}
(T\times U^{\lambda}_0)=U^{\lambda}$,
such that 
$\rm R^i(\bar{\mu}|_{U^{\lambda}})_{\ast}(W|_{U^{\lambda}}(D^{\lambda}))=0$
for all $\rm i>0$ and all $\rm \lambda$.
We can choose $\rm D^{\lambda}$ of the form $\rm T\times D_0^{\lambda}$,
where $\rm D_0^{\lambda}$ is a sufficiently high multiple of a section over
$\rm U_0^{\lambda}$ meeting the fibres of $\rm \mu$ transversally.
Since M is compact we may assume that 
$\rm \{ U^{\lambda}_0|\, \lambda =1,\ldots ,n\}$ already covers M.
The short exact sequence
\[ \rm 0\rightarrow W|_{U^{\lambda}} \rightarrow W|_{U^{\lambda}}(D^{\lambda}) \rightarrow W|_{U^{\lambda}}(D^{\lambda})\otimes {\cal O}_{D^{\lambda}} \rightarrow 0 \]
yields the exact sequence
\[ \rm 0\rightarrow (\bar{\mu}|_{U^{\lambda}})_{\ast}W|_{U^{\lambda}} 
\rightarrow (\bar{\mu}|_{U^{\lambda}})_{\ast}W|_{U^{\lambda}}(D^{\lambda}) \stackrel{u^{\lambda}}{\longrightarrow} 
(\bar{\mu}|_{U^{\lambda}})_{\ast}W|_{U^{\lambda}}(D^{\lambda})
\otimes {\cal O}_{D^{\lambda}} \rightarrow R^1(\bar{\mu}|_{U^{\lambda}})_{\ast}W|_{U^{\lambda}} \rightarrow 0 .\]

In the following h denotes a point on H as well as the corresponding line
on P. If $\rm V_t$ is trivial along h, then $\rm h\cong \PP^1_{\CC}$ gives
$\rm W|_{\bar{\mu}^{-1}(t,h)}\cong {\cal O}_{\PP^1_{\CC}}(-1)$, since we had
assumed $\rm (h\cdot S)=1$. Since $\rm V_0$ is trivial along twistor fibres
and since the triviality of a vector bundle on $\rm \PP^1_{\CC}$ is an open
property, there is an open neighbourhood 
$\rm \{ 0 \} \times M \subset B \subset T\times H$, such that for all
$\rm (t,h)\in B$ the vector bundle 
$\rm V_t|_h$ is trivial.
Inparticular we have  for all $\rm (t,h)\in B$
\[ \rm H^0(\bar{\mu}^{-1}(t,h),W|_{\bar{\mu}^{-1}(t,h)})=
       H^1(\bar{\mu}^{-1}(t,h),W|_{\bar{\mu}^{-1}(t,h)})=0. \]

Since $\rm \bar{\mu}$ is flat, W is flat over $\rm T\times H$
and base change implies the vanishing of 
$\rm (\bar{\mu}|_{U^{\lambda}})_{\ast}W|_{U^{\lambda}}$ and
of $\rm R^1(\bar{\mu}|_{U^{\lambda}})_{\ast}W|_{U^{\lambda}}$ over 
$\rm B\cap (T\times U^{\lambda}_0)$.
Therefore $\rm u^{\lambda}$ is an isomorphism over 
$\rm B\cap (T\times U^{\lambda}_0)$
and hence $\rm det(u^{\lambda})$ is a not constantly zero section of
\[ \rm 
det\! \left(
((\bar{\mu}|_{U^{\lambda}})_{\ast}\!W|_{U^{\lambda}}(D^{\lambda}))^{\vee}\!
\otimes
((\bar{\mu}|_{U^{\lambda}})_{\ast}\!W|_{U^{\lambda}}(D^{\lambda})
\otimes {\cal O}_{D^{\lambda}})
\right)
\cong
det\! \left(
((\bar{\mu}|_{U^{\lambda}})_{\ast}\!W|_{U^{\lambda}})^{\vee}\!
\otimes
R^1(\bar{\mu}|_{U^{\lambda}})_{\ast}\!W|_{U^{\lambda}}
\right).
\]
Moreover, all points $\rm (t,h)$ with $\rm V_t|_h$ non-trivial
belong to the support of $\rm \Delta^{\lambda}=div(det\, u^{\lambda})$.

For $\rm \lambda = 1,\ldots ,n\,$ there are open neighbourhoods
$\rm 0\in T_0^{\lambda}\subset T$, such that 
$\rm \Delta^{\lambda}_t = \Delta^{\lambda}|_{\{ t\} \times U_0^{\lambda}}$
are divisors on $\rm U_0^{\lambda}$ for all $\rm t\in T^{\lambda}_0$.
Thus $\rm \Delta^{\lambda}\subset T^{\lambda}_0 \times U_0^{\lambda}$
are analytic families of divisors on $\rm U_0^{\lambda}$ for all
$\rm \lambda = 1,\ldots ,n$.
Since $\rm V_0$ is an instanton bundle we have
$\rm \Delta^{\lambda}_0\subset U_0^{\lambda}-M$. 
M is a closed subset of H and thus there exists
open neighbourhoods $\rm 0\in T_1^{\lambda}\subset T_0^{\lambda}$,
such that for all $\rm t\in T^{\lambda}_1$ we have 
$\rm \Delta^{\lambda}_t\subset U_0^{\lambda}-M$.
We put $\rm T_1 = \bigcap_{\lambda = 1}^{n}T_1^{\lambda}$ and obtain
that for all $\rm t\in T_1$ and for all $\rm h\in M$ the vector bundle
$\rm V_t|_h$ is trivial. $\Box$

Warning: Contrarily to the most other parts of this paper,
the last lemma does not hold in the algebraic setup, since M has not
to be Zariski-closed.

\noindent
{\it Remark.}
Let $\rm {\cal M}_{\cal I}(P,F,{\cal O}_F^r)$ be the functor of families of
framed mathematical instanton bundles.
By the previous Lemma we have have
$\rm {\cal M}_{\cal I}(P,F,{\cal O}_F^r) \subset {\cal M}(P,F,{\cal O}_F^r)$
as an analytically open subfunctor.
In particular,
\[ \rm {\cal M}_1(P,F,{\cal O}_F^r,\pi^{\ast}c_{\bullet})
       ={\cal M}_{\cal I}(P,F,{\cal O}_F^r)
       \cap {\cal M}_0(P,F,{\cal O}_F^r,\pi^{\ast}c_{\bullet}) \]
is an analytically open subfunctor and the natural transformation
\[ \rm {\cal M}_1(P,F,{\cal O}_F^r,\pi^{\ast}c_{\bullet}) \rightarrow
       {\cal M} (S,F,{\cal O}_F^r,\pi_S^{\ast}c_{\bullet}) \]
is smooth of relative dimension
$\rm -\chi(\End (V)(-S))$ with all the details as in Theorem 2.9.
\indent

For T an algebraic space and 
$\rm (V,\alpha)\in {\cal M}_{\cal I}(P,F,{\cal O}_F^r)(T)$
we define 
$\rm \sigma(V,\alpha)$
to be the family of framed instanton bundles
$\rm ((id_T\times \tau)^{\ast}\bar{V}^{\vee},
            (\tau^{\ast}\bar{\alpha}^{\vee})^{-1})$,
where $\tau$ is the antiholomorphic involution on P as in the introduction
and the "bar" means complex conjugation.
Thus $\sigma$ defines a real structure on
$\rm {\cal M}_{\cal I}(P,F,{\cal O}_F^r)$ and the closed fixpoints
$\rm {\cal M}_{\cal I}(P,F,{\cal O}_F^r)(\RR)$ is just the set of isomorphism classes of framed $\rm U(r)-instanton$ bundles on P.

To be more specific, we define $\rm {\cal M}_{\cal I}^{\RR}(P,F,{\cal O}_F^r)$
as the functor that associates to an algebraic space T
the set of isomorphism classes of families of framed instanton bundles
$\rm (V,\alpha)$ on P parametrized by T, which have the additional property
that there exists an isomorphism of framed vector bundles
\[ \rm
\varphi : (V,\alpha) \rightarrow ((id_T\times \tau)^{\ast}\bar{V}^{\vee},
            (\tau^{\ast}\bar{\alpha}^{\vee})^{-1}) .\]
Then $\sigma$ defines 
$\rm {\cal M}_{\cal I}^{\RR}(P,F,{\cal O}_F^r)$
as closed real-analytic subfunctor of $\rm {\cal M}_{\cal I}(P,F,{\cal O}_F^r)$.

We have $\rm {\cal M}_{\cal I}^{\RR}(P,F,{\cal O}_F^r)$
embedded in $\rm {\cal M}_1(P,F,{\cal O}_F^r)$
and we can examine the restriction of the natural transformation
described in Section 2.2 to 
$\rm {\cal M}_{\cal I}^{\RR}(P,F,{\cal O}_F^r)$.
\begin{theorem}
The natural transformation
\[ \rm {\cal M}_{\cal I}^{\RR}(P,F,{\cal O}_F^r) \rightarrow
{\cal M}(S,F,{\cal O}_F^r) \]
given by restriction is an open real-analytic embedding.
\end{theorem}

{\it Proof:} We have to show that the restriction is injective
for closed points and the associated tangent maps are isomorphisms.

We denote $\rm \tau |_{S}$ and 
$\rm \tau |_{\bar{S}}$ again
with $\tau$ and obtain a real structure $\sigma'$ on
$\rm {\cal M}(S,F,{\cal O}_F^r)\times {\cal M}(\bar{S},F,{\cal O}_F^r)$
by mapping 
\[ \rm \left( (V,\alpha),(W,\beta) \right) \in 
{\cal M}(S,F,{\cal O}_F^r)(T)\times {\cal M}(\bar{S},F,{\cal O}_F^r)(T)\]
to
\[ \rm \left( ((id_T\times \tau)^{\ast}\bar{W}^{\vee},
            (\tau^{\ast}\bar{\beta}^{\vee})^{-1}),
((id_T\times \tau)^{\ast}\bar{V}^{\vee},
            (\tau^{\ast}\bar{\alpha}^{\vee})^{-1}) \right).\]
Hence the natural transformation
\[ \rm \Psi : {\cal M}_{\cal I}(P,F,{\cal O}_F^r) \rightarrow
{\cal M}(S,F,{\cal O}_F^r)\times {\cal M}(\bar{S},F,{\cal O}_F^r) \]
given by restriction to both factors is real, i.e. compatible with both real structures $\sigma$ and $\sigma'$. 

From the diagram $\rm (\ast)$ in the proof of Theorem 2.9 we obtained
for $\rm (V,\alpha) \in {\cal M}_0 (P,F,{\cal O}_F^r)(\CC)$
the commutative diagram
\[ \begin{array}{ccc}
0&&0\\
\downarrow & & \downarrow \\
\rm H^1(P,\End (V)(-S))&\stackrel{=}{\longrightarrow}&\rm H^1(P, \End (V)(-S))\\
\downarrow & & \downarrow \\
\begin{array}{c}
 \rm H^1(P,\End (V)(-S))\\[-1ex]
 \oplus\\[-1ex]
 \rm H^1(P,\End (V)(-\bar{S}))
\end{array}&
\stackrel{\cong}{\longrightarrow}&\rm H^1(P, \End (V)\otimes {\cal J}_{F|P})\\
\downarrow & & \downarrow \\
\rm H^1(P,\End (V)(-\bar{S}))&\stackrel{\cong}{\longrightarrow}&
 \rm H^1(S,\End (V|_S)(-F))\\
\downarrow & & \downarrow \\
0&&0
\end{array} \]
with exact columns.
We may interchange $\rm S$ and $\rm \bar{S}$ and make use of
the isomorphism
\[ \rm H^1(P,\End (V)(-S)) \stackrel{\cong}{\longrightarrow} 
 H^1(\bar{S},\End (V|_{\bar{S}})(-F)) \]
to obtain a natural splitting
\[ \rm H^1(P, \End (V)\otimes {\cal J}_{F|P}) =
H^1(S,\End (V|_S)(-F)) \oplus H^1(\bar{S},\End (V|_{\bar{S}})(-F)), \]
which is just the tangent map of $\rm \Psi$ at $\rm (V,\alpha)$.

Now we assume that 
$\rm (V,\alpha)$ and $\rm (W,\beta)$ are two elements in
$\rm {\cal M}_{\cal I}(P,F,{\cal O}_F^r)(\RR)$, such that
$\rm \Psi(V,\alpha)\cong \Psi(W,\beta)$. This means that there are two
isomorphisms of framed vector bundles
$\rm \varphi_S : (V|_S,\alpha)\rightarrow (W|_S,\beta)$ and
$\rm \varphi_{\bar{S}} : (V|_{\bar{S}},\alpha)\rightarrow (W|_{\bar{S}},\beta)$, which implies that 
$\rm \varphi_S|_F = \beta^{-1}\circ \alpha = \varphi_{\bar{S}}|_F$.
Thus we have an isomorphism of framed vector bundles
$\rm \varphi_{S\cup \bar{S}} : (V|_{S\cup \bar{S}},\alpha)\rightarrow (W|_{S\cup \bar{S}},\beta)$. With V and W U(r)-instanton bundles also
$\rm \Hom (V,W)$ is an U(r)-instanton bundle and in particular
$\rm H^1(P,\Hom (V,W)(-S-\bar{S}))=0$. Therefore the restriction map
$\rm Hom(V,W)\rightarrow Hom(V|_{S\cup \bar{S}},W|_{S\cup \bar{S}})$ is
surjective and we can find a morphism of framed vector bundles
$\rm \varphi : (V,\alpha) \rightarrow (W,\beta)$ as extension of
$\rm \varphi_{S\cup \bar{S}}$. By Lemma 2.3 $\rm \varphi$ is an isomorphism
and hence $\Psi$ injective along 
$\rm {\cal M}_{\cal I}(P,F,{\cal O}^r_F)(\RR)$.

Therefore $\rm \Psi$ is an open embedding in an open neighbourhood
of $\rm {\cal M}_{\cal I}^{\RR}(P,F,{\cal O}^r_F)$ and hence
we have a real-analytic open embedding
\[ \rm {\cal M}_{\cal I}^{\RR}(P,F,{\cal O}^r_F) \rightarrow {\cal N}^{\RR}, \]
where $\rm {\cal N}^{\RR}$ denotes the closed real-analytic subfunctor
of $\rm {\cal M}(S,F,{\cal O}_F^r)\times {\cal M}(\bar{S},F,{\cal O}_F^r)$
defined by the real structure $\rm \sigma'$.

Now we consider the projection 
\[ \Phi :\rm {\cal M}(S,F,{\cal O}_F^r)\times {\cal M}(\bar{S},F,{\cal O}_F^r)
\rightarrow {\cal M}(S,F,{\cal O}_F^r) \]
and notice that $\rm \Phi $ induces an isomorphism
$\rm {\cal N}^{\RR} \rightarrow {\cal M}(S,F,{\cal O}_F^r)$ due to the inverse
transformation
\[ \begin{array}{ccc}
\rm {\cal M}(S,F,{\cal O}_F^r)(T) & \longrightarrow &
\rm {\cal N}^{\RR}(T)\hookrightarrow \left( {\cal M}(S,F,{\cal O}_F^r)\times {\cal M}(\bar{S},F,{\cal O}_F^r) \right) (T)\\
\rm (V,\alpha)& \longrightarrow & \rm \left( (V,\alpha), 
((id_T\times \tau)^{\ast}\bar{V}^{\vee},
            (\tau^{\ast}\bar{\alpha}^{\vee})^{-1})  \right).
\end{array} \]
Thus $\rm \Phi \circ \Psi$ defines an open real-analytic embedding
\[ \rm {\cal M}_{\cal I}^{\RR}(P,F,{\cal O}_F^r) \rightarrow
{\cal M}(S,F,{\cal O}_F^r) .\;\;\Box \]

The result of N.P. Buchdahl's article \cite{buchdahl2} implies that
$\rm {\cal M}_{\cal I}(P,F,{\cal O}_F^r)(\RR)$ and
$\rm {\cal M}(S,F,{\cal O}_F^r)(\CC)$ are naturally bijective and
we obtain 
\begin{corollary}
The natural transformation
$ \rm {\cal M}_{\cal I}^{\RR}(P,F,{\cal O}_F^r) \rightarrow
{\cal M}(S,F,{\cal O}_F^r) $
given by restriction is a real-analytic isomorphism.
For fixed $\rm c_{\bullet}\in H^{\bullet}(M,\ZZ)$ the functor
$\rm {\cal M}_{\cal I}^{\RR}(P,F,{\cal O}_F^r,\pi^{\ast}c_{\bullet})$
is represented by a smooth separated algebraic space of finite type
over $\CC$.
\end{corollary}

\renewcommand{\baselinestretch}{1.0}

\end{document}